\documentclass[12pt]{article}

\usepackage{a4wide}
\usepackage{amsmath}
\usepackage{amsfonts}
\usepackage{amssymb}
\usepackage{latexsym}
\usepackage{verbatim}
\usepackage{curves}

\newtheorem{theorem}{Theorem}

\newtheorem{remark}[theorem]{Remark}

\newtheorem{proposition}[theorem]{Proposition}

\newcommand{\R}{\hbox{\rm I \kern-.50em R}}
\newcommand{\Z}{\hbox{\rm Z \kern-.65em Z}}

\begin{document}

\baselineskip15pt


\title{On Conjugacy Classes of Groups of Squarefree Order}
{\medskip}

\author{
Anna Torstensson \\
{\normalsize Centre for Mathematical Sciences} \\
{\normalsize Box 118, SE-221 00 Lund, Sweden} \\
{\normalsize Email: anna.torstensson@math.lth.se}}

\date{}
\maketitle
 \begin{abstract}
    The problem of finding the largest finite group with a certain class
    number (number of conjugacy classes), $k(G)$, has been investigated by a number
of researchers since the early 1900's and has been solved by computer for
$k(G) \leq 9$. (For the restriction to simple
groups for $k(G) \leq 12$.) One has also tried to find a general upper bound on $|G|$ in
terms of $k(G)$. The best known upper bound in the general case is
in the order of magnitude $|G| \leq k(G)^{2^{k(G)-1}}$.
In this paper we consider the restriction of this longstanding problem to
groups of squarefree order. We derive an explicit formula for the class number $k(G)$ of any
group of squarefree order and we also obtain an estimate $|G| \leq
k(G)^3$ in this case. Combining the two results we get an
efficient way to compute the largest squarefree order a group with a
certain class number can have. We also provide an implementation of this
algorithm to compute the maximal squarefree $|G|$ for arbitrary $k(G)$
and display the results obtained for $k(G) \leq 100$.

\end{abstract}


%
%
\section{Introduction}

Conjugacy classes are a vital tool in the analysis of groups. A
condensed piece of information about a group is its class number, that
is its number of conjugacy classes. Starting with small numbers
mathematicians have tried to catalogise all groups having a given
class number. A very brief history of this problem with references can
be found in \cite{Poland}. The starting point was when Landau in 1903
managed to show that there is only a finite number of
non-isomorphic finite groups with class number $k$, and also derived an
upper bound on the group size as a function of $k$. In \cite{Poland}
Poland extends the complete list of groups with class number
$k$ from $k \leq 5$ to $k \leq 7$. Using computerised methods this was
later pushed forward to cover $k \leq 9$ (references can be found in
\cite{Komissarcik}) and considering only simple groups further to $k
\leq 12$ by Komissarcik in \cite{Komissarcik}. The reason to limit the
search to simple groups was that for $5 \leq k \leq 9$
the largest groups with a given class number are all simple. Hence one
might guess that this is true also for larger $k$ even though this
seems to be unknown.

There has also been successful attempts to attack
other types of groups. For example $p$-groups where it is known that a
group of order $p^m$ has $(n+r(p-1))(p^2-1)+p^e$ conjugacy classes
where $r$ is some non-negative integer and $m=2n+e$, $e= 0$ or $1$ depending on the parity
of $m$. This result can be found in \cite{Poland2}.

In this paper we
will go in the opposite direction and instead of $p$-groups consider
groups of squarefree order. We will derive an explicit formula for the
class numbers groups of a given squarefree order can have. As a
by product of the analysis we will also see that a group of squarefree
order having $k$ conjugacy classes can have no more than $k^3$ elements. As
mentioned there are bounds on the group order in the general case, but
in contrast to this bound they are far from sharp.

It is of course a substantial confinement to consider only groups of
squarefree orders. It should be noted however that most group orders
are covered, as $\lim_{n \to \infty} \frac{Q(n)}{n} =
\frac{6}{\pi^2}$, where $Q(n)$ denotes the number of squarefree
numbers $\leq n$. (See \cite{HardyWright}.) In addition the techniques
developed here can hopefully be modified in order to cover other
types of groups which have relatively simple presentations.

\section{Description of the conjugacy classes}

In this section we will describe the conjugacy classes of groups of
squarefree order. Along the way we will also introduce some notation
that will be used throughout this article.

Let $G$ be a group of squarefree order $n$. Then Murty and Murty
(\cite{MurtyMurty}) have proven that $G$ has a
presentation of the form

\begin{equation}\label{Gpres} \langle x,y|x^d=id,y^e=id,x^{-1}yx=y^t \rangle \end{equation}

for some divisor $d$ of $n$, $e=n/d$ and some integer $t$ which is of order
$d$ modulo $e$. We will use this presentation to compute the conjugacy
classes of $G$.

\begin{remark}\label{decond} Note that an element $t$ of order $d$ modulo $e$ exists
  only for those factorisations $n=d \cdot e$ which satisfy
  $(d,\phi(e))=d$. This is since $(d,\phi(e))=d'<d$ implies
  $t^{d'} \equiv t^{ad+b\phi(e)} \equiv 1$ contradicting that $t$ is
  order $d$.
\end{remark}

\begin{proposition}\label{classsize} Let $G$ be a group of squarefree order $n=de$ and
$$ \langle x,y|x^d,y^e,x^{-1}yx=y^t \rangle$$
a presentation of $G$. Then the size of the conjugacy
class containing the element $x^ay^b$ is $k{e \over
  (e,1-t^a)}$ where $k$ is the smallest positive integer for which $b(t^k-1) \equiv 0$ modulo $(e,1-t^a)$.
\end{proposition}

{\bf Proof:}
Any element of $G$ has a unique representation
$x^ay^b$ with $0 \leq a \leq d-1$ and $0 \leq b \leq e-1$ and a
straightforward computation using the relations gives that

$$ x^{-1}x^ay^bx=x^ay^{tb} \hspace{7mm} \mbox{ and } \hspace{7mm} y^{-1}x^ay^by=x^ay^{b+1-t^a} $$

This shows that elements of the form $x^ay^b$ can be conjugates only
of elements of the form $x^ay^{b_1}$. Moreover these two elements are
conjugate if and only if $b_1$ can be obtained from $b$ by a sequence
of the operations $O_1: b \mapsto tb$ and $O_2: b \mapsto b+1-t^a$. Let $S_b$ be the
set obtained from the element $b$ in this way. The size of $S_b$ equals the size of the
conjugacy class containing $x^ay^b$ so we would like to compute
$|S_b|$. In fact we can describe the set $S_b$ in the following way:
Let k be the smallest positive integer such that

$$b(t^k-1) \equiv s(1-t^a) \mbox{ mod } e $$

for some integer $s$. Note that such a $k$ exists and is at most $a$
since for $k=a$ we may take $s=-b$. Now we will show that

$$ S_b=\{t^rb+s(1-t^a)| 0 \leq r \leq k-1, s \in \mathbb{Z} \} $$

We regard $S_b$ as a subset of $\mathbb{Z}_e$ as it is only the
residue modulo $e$ that affects which element $x^ay^b$ we get in
$G$. It is clear that an arbitrary element $t^rb+s(1-t^a)$ in

$$ M_b=\{t^rb+s(1-t^a)| 0 \leq r \leq k-1, s \in \mathbb{Z} \} $$

can be obtained from $b$ simply by first applying $b \mapsto tb$ $r$
times and then $b \mapsto b+1-t^a$ $s$ times. Thus it is clear that $M_b
\subseteq S_b$. On the other hand $M_b$ contains $b$ and it is
straightforward to check that it is closed under the two operations
$O_1$ and $O_2$. Hence we must have $M_b=S_b$. Also two elements
in
$$ S_b=\{t^rb+s(1-t^a)| 0 \leq r \leq k-1, s \in \mathbb{Z} \}
\subseteq \mathbb{Z}_e $$

are different $t^{r_1}b+s_1(1-t^a) \not \equiv
t^rb+s(1-t^a)$ mod $e$ exactly when $s \not \equiv s_1$ modulo ${e \over (e,1-t^a)}$ since if we have
$$t^rb+s(1-t^a) \equiv t^{r_1}b+s_1(1-t^a) \mbox{ mod } e$$ with $0 \leq r_1
\leq r \leq k-1$ then $t^r_1(b-t^{r-r_1}b) \equiv (s-s_1)(1-t^a)$ and
as $t$ is order $d$ mod $e$ $t$ is invertible and we get $(b-t^{r-r_1}b)
\equiv t_1^{-r}(s-s_1)(1-t^a)$. From the definition of $k$ it follows
that $r-r_1 \geq k$ unless $r=r_1$ so considering the condition $0 \leq r_1
\leq r \leq k-1$ the numbers $r$ and $r_1$ must be equal. It follows that $s(1-t^a) \equiv
s_1(1-t^a)$ modulo $e$ and this holds if and only if $s \equiv s_1$ modulo ${e \over
  (e,1-t^a)}$. This shows that there are $k{e \over
  (e,1-t^a)}$ different elements in the set $S_b$.

The proposition now follows by
noting that there is an integer $s$ with
$$b(t^k-1) \equiv s(1-t^a) \mbox{ mod } e $$
if and only if $b(t^k-1) \equiv 0$ modulo $(e,1-t^a)$.

$\blacklozenge$

\vspace{4mm}

In the following we shall use the smallest common divisors $(e,1-t^s)$
frequently and therefore we introduce the notation $\alpha_s=(e,1-t^s)$
for $1 \leq s \leq d$.

As a consequence of the above proposition we get a lower bound on the
number of conjugacy classes of a group of squarefree order.

\begin{theorem}\label{upperbound}  Let $G$ be a group of squarefree order and $k(G)$ the
  number of conjugacy classes of $G$. Then $|G| \leq k(G)^3$.
\end{theorem}

{\bf Proof:} Let $n=de$ be any factorisation of $n=|G|$. By the above
proposition we know that any conjugacy class containing an element
$x^ay^b$ has $k{e \over \alpha_a}$ elements for some $k$ between
$1$ and $d$. Thus, if we denote the class containing $x$ by $C(x)$,  ${e \over \alpha_a} \leq |C(x^ay^b)| \leq d{e \over
  \alpha_a} $. For a fixed $a$ there are $e$ elements to distribute
among such classes and hence the resulting number of classes must be
between ${\alpha_a \over d}$ and $\alpha_a$. Summing over all $a$ we
get ${1 \over d}\sum_{a=1}^{d}\alpha_a \leq k(G) \leq
\sum_{a=1}^{d}\alpha_a$. Noting that $\alpha_d=(e,1-t^d)=e$ we get
${e \over d} \leq {1 \over d}\sum_{a=1}^{d}\alpha_a \leq k(G)$. On
the other hand each $a$ results in at least one class so $d \leq
k(G)$. Taken together $|G|=n=de={e \over d}d^2 \leq k(G)^3$.

$\blacklozenge$

\vspace{4mm}

\begin{remark} This bound does not hold for groups of arbitrary
orders. In fact a smallest counter example is given by $A_6$, the
alternating group on 6 symbols, which has $360$ elements but only $7$
conjugacy classes. \end{remark}

This bound can be compared with the smallest bound known to hold for
all groups which is in the order of magnitude $|G| \leq k(G)^{2^{k(G)-1}}$. (See
\cite{Newman}.) There are however lower estimates for a number of
specific types of groups.

\section{A formula for the class number}\label{formulasec}

Let us continue to compute the number of conjugacy classes for groups
of squarefree order, which as we have seen always have some
presentation of the form (\ref{Gpres}). Let us for the rest of this section assume that $d$, $e$ and $t$ are fixed,
in other words that we are working with a fixed group and a certain
presentation of it. We also assume that $n$ is squarefree and hence that $e$ can be written as a
product $e=p_1p_2p_3 \cdots p_m$ of distinct primes. From Proposition
\ref{classsize} it is clear that the sizes of the conjugacy classes
depend only on the sequence ${\alpha} = (\alpha_1, \alpha_2, \ldots ,
\alpha_d)$. On the other hand $\alpha_s=(e,1-t^s)$ is a product of
those primes $p_i$ for which the order of $t$ modulo $p_i$ divides
$s$. Hence $\alpha$ can be directly computed from the order of $t$
modulo the divisors $p_i$ of $e$. The fact that $t$ is of order $d$
modulo $e$ restricts the possible orders modulo $p_i$ to numbers
dividing $(d,p_i-1)$, but for any $d_i|(d,p_i-1)$ there is a $t$ with
this order modulo $p_i$. From the Chinese Remainder Theorem it follows
that for any sequence ${\bf d}=(d_1, d_2, \ldots , d_r)$ with  $d_i|(d,p_i-1)$
there is a $t$ such that $t$ has order $d_i$ modulo $p_i$ for each
$i$. From the sequence ${\bf d}$ we now get the corresponding ${\bf \alpha}$-sequence
by $\alpha_l=\prod_{d_i|l}p_i$. Let us for future use note two
interesting properties of ${\bf \alpha}$ that follows from this description:

\begin{proposition}\label{alphagcd} Let $e$ be a squarefree number and $t$ a
  positive integer. The sequence  ${\bf \alpha} = (\alpha_1, \alpha_2, \ldots ,
\alpha_d)$ where $\alpha_s=(e,1-t^s)$ then has the property that
$(\alpha_k,\alpha_l)=\alpha_{(k,l)}$.
\end{proposition}

{\bf Proof:} As above we let $d_i$ be the order of $t$ modulo $p_i$ and
thus get the expression  $\alpha_l=\prod_{d_i|l}p_i$ for
$\alpha_l$. It follows that
$$
(\alpha_k,\alpha_l)=(\prod_{d_i|k}p_i,\prod_{d_i|l}p_i)=\prod_{d_i|(k,l)}p_i
= \alpha_{(k,l)}$$

$\blacklozenge$

\vspace{4mm}

\begin{proposition}\label{alphagcdwithd} Let $e$ be a squarefree number and $t$ some positive integer which is of order $d$ modulo $e$. The sequence  ${\bf \alpha} = (\alpha_1, \alpha_2, \ldots ,
\alpha_d)$ where $\alpha_s=(e,1-t^s)$ then has the property that
$\alpha_k=\alpha_{(k,d)}$.
\end{proposition}

{\bf Proof:}  As before we know that $\alpha_l=\prod_{d_i|l}p_i$ where
$d_i$ is the order of $t$ modulo $e$. Clearly we must have that
$d_i|d$ since $d$ is the order of $t$ modulo $e=p_1p_2p_3 \cdots
p_m$. It follows that $d_i|l$ if and only if $d_i|(l,d)$. Consequently

$$\alpha_l=\prod_{d_i|l}p_i=\prod_{d_i|(d,l)}p_i=\alpha_{(d,l)}$$

$\blacklozenge$

\vspace{4mm}

Above we have constructed the
sequence ${\bf \alpha}$ for a given group of squarefree order $n=de$. Now we want to compute the number of
conjugacy classes from the sequence ${\bf \alpha}$. By Proposition
\ref{classsize} we know that the size of the
class containing $x^ay^b$ is $k{e \over \alpha_a}$ and we observe that
k is the smallest positive integer such that $b \equiv 0$ modulo
${\alpha_a \over (\alpha_a,\alpha_k)}$. To count how many $b$ we get
for a certain value of $k$ we introduce the following notation. Let
$$C_j=\{b \hspace{2mm}| \hspace{2mm} b \equiv 0 \hspace{4mm} \mbox{modulo} \hspace{4mm}
{\alpha_a \over (\alpha_a,\alpha_j)}\} \subseteq \mathbb{Z}_e$$
The elements $x^ay^b$ ending up in conjugacy classes of size $k{e
  \over \alpha_a}$ are those where $b$ is in the set $C_k \setminus
\cup_{j=1}^{k-1} C_j$. From the principle of inclusion and exclusion
we get an expression for the size of this set that can be used for computation.
$$|C_k \setminus \bigcup_{j=1}^{k-1}C_j|=
  |C_k|+\sum_{l=1}^{k-1}(-1)^l \sum_{j_1<j_2< \cdots <j_l<k}|C_k \cap
  \bigcap_{s=1}^l C_{j_s}|$$

Now from the definition of $C_k$ and proposition \ref{alphagcd} it
follows that

$$|\bigcap_{s=1}^t
C_{j_s}|=\frac{e(\alpha_a,\alpha_{j_1},\alpha_{j_2}, \ldots,
  \alpha_{j_t})}{\alpha_a}=\frac{e\alpha_{(a,j_1,j_2 \ldots j_t)}}{\alpha_a}$$

and hence the above expression can be rewritten as

$$|C_k \setminus \bigcup_{j=1}^{k-1}C_j|=
  \frac{e}{\alpha_a}(\sum_{l=0}^{k-1}(-1)^l
  \sum_{j_1<j_2< \cdots <j_l<k} \alpha_{(a,k,j_1,j_2, \ldots ,
    j_l)})$$
Now the above expression counts the number of elements in the group
containing exactly $a$ factors of $x$ that belong to a conjugacy class
of size $k{e \over \alpha_a}$. We can now compute the total number of
classes for fixed $k$:

$$  \sum_{a=1}^d(\frac{1}{k}(\sum_{l=0}^{k-1}(-1)^l
  \sum_{j_1<j_2< \cdots <j_l<k} \alpha_{(a,k,j_1,j_2, \ldots ,
    j_l)}))$$

Finally summation over all possible values of $k$ results in a formula
for the total number of conjugacy classes in the group:

$$  k(G)=\sum_{k=1}^d\sum_{a=1}^d(\frac{1}{k}(\sum_{l=0}^{k-1}(-1)^l
  \sum_{j_1<j_2< \cdots <j_l<k} \alpha_{(a,k,j_1,j_2, \ldots ,
    j_l)}))$$

The above formula is explicit but not very handy for computation. We
therefore rewrite it into a form where all $\alpha_j$ are collected
together.

\begin{theorem}
  Let $n=de$ be a squarefree number with $e=p_1p_2p_3 \cdots p_m$ and $G$
  the group with presentation (\ref{Gpres}). Further, let $d_i$ be the
  order of $t$ modulo $p_i$ and $\alpha_j=\prod_{d_i|j}p_i$. Then the
  class number of $G$ equals

  $$k(G)=\sum_{k=1}^d\sum_{a=1}^d(\frac{1}{k}(\sum_{l=0}^{k-1}(-1)^l
  \sum_{j_1<j_2< \cdots <j_l<k} \alpha_{(d,a,k,j_1,j_2, \ldots ,
    j_l)}) = $$
  \begin{equation}\label{formula}= d
    \sum_{j|d}\frac{\alpha_j}{j^2}\sum_{s|\frac{d}{j}}
    \frac{\mu(s)}{s^2} \end{equation}
  where $\mu$ denotes the M{\"o}bius function.
\end{theorem}

\begin{remark} The first expression for $k(G)$ is clear from the
discussion above if we only note that it is legitimate to include $d$
in the greatest common divisors by proposition
\ref{alphagcdwithd}. The reason for doing that is to simplify the
formula by reducing the number of terms leaving only those $\alpha$ indices which are
divisors of $d$.
\end{remark}

{\bf Proof:} The proof is by induction on the indices of
$\alpha$. First we prove that for any $d$ the coefficient of $\alpha_1$
is equal on both sides of the equation.

Let $d=q_1q_2 \ldots q_t$. First look at the inner sum of the LHS for
a fixed pair $(a,k)$. Let $(a,k,d)=q_1q_2 \ldots q_s$ (where $s \leq
t$). We are interested in the coefficient of $\alpha_1$ so we
should, for each $l$, count the number of sequences $(j_1,j_2, \ldots, j_l)$ with $(a,k,d,j_1,j_2, \ldots, j_l)=1$. Now  $(a,k,d,j_1,j_2, \ldots, j_l)=1$ if and only if
none of the primes $(a,k,d)=q_1q_2 \ldots q_s$ divide every number in
the sequence $(j_1,j_2, \ldots, j_l)$. Using the principle of
inclusion and exclusion we can count the number of such sequences. Let
$D_i=\{(j_1,j_2, \ldots, j_l)\hspace{2mm}|\hspace{2mm} q_i|j_r \mbox{ for }r=1,2, 3 \ldots, l \}$,
the set of sequences divisible by $q_i$ and $U=\{(j_1,j_2, \ldots,
j_l)\}$, the set of all sequences. Then

$$|U \backslash \bigcup D_i|=|U|+\sum_{m=1}^s (-1)^m \sum_{1 \leq i_1<i_2< \ldots
  < i_m \leq s} |\bigcap_{r=1}^m D_{i_r}|$$
Noting that

$$|\bigcap_{r=1}^m D_{i_r}|=\left( \begin{array}{c}
    \frac{k}{q_{i_1}q_{i_2} \cdots q_{i_m}} -1 \\ l \\
  \end{array} \right)$$

as the elements of these sequences are chosen among the
$\frac{k}{q_{i_1}q_{i_2} \cdots q_{i_m}}-1$ numbers between $1$ and $k-1$
that are divisible by $q_{i_1}q_{i_2} \cdots q_{i_m}$, we obtain

$$ |U \backslash \bigcup D_i|=\left( \begin{array}{c}
   k -1 \\ l \\
  \end{array} \right)+\sum_{m=1}^s (-1)^m \sum_{1 \leq i_1<i_2< \ldots
  < i_m \leq s} \left( \begin{array}{c}
    \frac{k}{q_{i_1}q_{i_2} \cdots q_{i_m}} -1 \\ l \\
  \end{array} \right) $$
We can now compute the coefficient of $\alpha_1$ in the inner sum of
the LHS of the proposition by summing over $l$ thus obtaining

$$ \sum_{l=0}^{k-1} (-1)^l \left( \left( \begin{array}{c}
   k -1 \\ l \\
  \end{array} \right)+\sum_{m=1}^s (-1)^m \sum_{1 \leq i_1<i_2< \ldots
  < i_m \leq s} \left( \begin{array}{c}
    \frac{k}{q_{i_1}q_{i_2} \cdots q_{i_m}} -1 \\ l \\
  \end{array} \right) \right) = $$
$$=\sum_{m=1}^s (-1)^m \sum_{1 \leq i_1<i_2< \ldots
  < i_m \leq s} \sum_{l=0}^{k-1}(-1)^l \left( \begin{array}{c}
    \frac{k}{q_{i_1}q_{i_2} \cdots q_{i_m}} -1 \\ l \\
  \end{array} \right)   $$

as the first sum vanishes by the binomial theorem. Here,
again by the binomial theorem, the innermost sum is zero unless
$k=q_{i_1}q_{i_2} \cdots q_{i_m}$ in which case the sum is 1. The
second case can occur only for $m=s$ as $(a,k,d)=q_1q_2 \ldots q_s$
and then the total sum is $(-1)^{s}$. Thus we have shown that the
coefficient of $\alpha_1$ in the inner
sum of the LHS equals $(-1)^{s}$ if $k=q_1q_2 \ldots q_s=(a,k,d)$ and zero
otherwise so we get the coefficient of $\alpha_1$ in the LHS as

$$\sum_{k=1}^d \sum_{a=1}^d \frac{1}{k}\left\{ \begin{array}{c}
    {(-1)^s}\hspace{3mm} \mbox{if} \hspace{3mm} k=q_1q_2 \ldots q_s \\ 0
    \hspace{5mm}\mbox{otherwise}\end{array} \right\} =$$
$$= \sum_{k|d} \frac{\mu(k)}{k}\frac{d}{k}$$
and this is just the coefficient of $\alpha_1$ in the RHS of the
proposition. This completes the proof of the base for the
induction.

Let us now move on to the induction step. We assume that the
coefficient of $\alpha_k$ in the LHS and the RHS of the proposition
are equal for each $k<j$ (for any $d$). We now want to show that the
coefficients of $\alpha_j$ are equal on both sides. Let $p$ be a prime
that divides $j$. We may assume that $p|d$ since otherwise the
coefficient of $\alpha_j$ is obviously zero on both sides of the equality.
Looking at the $LHS_d$ of the proposition we find that studying the
coefficient of $\alpha_j$ it is only interesting to sum over $a$ and
$k$ divisible by $p$. Writing $d=pd'$, $a=pa'$, $k=pk'$, $j=pj'$ and
$j_t=pj_t'$ we are looking for the coefficient of $\alpha_{pj'}$ in

$$ \sum_{a'=1}^{d'} \sum_{k'=1}^{d'} \frac{1}{pk'}
\sum_{l=0}^{pk'-1} (-1)^l \sum_{j_1'<j_2'< \cdots < j_l'<k'}
\alpha_{p(a',k',d',j_1',j_2', \ldots , j_l')}$$

Here the second last summation can just as well stop at $k'-1$ since
the last sum is empty for larger indices. Taking this into
consideration and looking at the above formula we find that the coefficient
of $\alpha_j=\alpha_{pj'}$ in $LHS_d$ equals $\frac{1}{p}$
times the coefficient of $\alpha_{j'}$ in $LHS_{d'}$. On the other
hand the latter equals $\frac{1}{p}$
times the coefficient of $\alpha_{j'}$ in $RHS_{d'}$ by the induction hypothesis. We would like to
show that this in turn equals the coefficient of $\alpha_{pj'}$ in $RHS_d$ so
let us compare those coefficients. The last one is

$$\frac{pd'}{(pj')^2} \sum_{s|\frac{d}{j}} \frac{\mu(s)}{s^2} =$$
$$=\frac{d'}{pj'^2} \sum_{s|\frac{d'}{j'}} \frac{\mu(s)}{s^2} $$

which equals just$\frac{1}{p}$ times the coefficient of $\alpha_{j'}$
in $RHS_{d'}$. This completes the induction step and thus the proof of
the theorem.

$\blacklozenge$

\section{Computation of class numbers}

Using the class number formula above one can easily compute all
possible class numbers for groups of a given squarefree order
$n$. The algorithm was implemented in the computer algebra system
Magma (see appendix).

The resulting class numbers for small group
orders $n$ are shown in Table \ref{tabellall}.

\noindent
\begin{table}
  \caption{All class numbers a group of given squarefree order can
    have. Prime orders have been omitted to save space as the only possible class
    number in this case is the group order.} \label{tabellall}
   \vspace{2mm}
\begin{tabular}{|l|l|l|l|l|l|}
  \hline
  $|G|$ & Smallest $k(G)$ & All $k(G)$ & $|G|$ & Smallest $k(G)$ & All $k(G)$ \\
  \hline \hline
  1 & 1 & 1 & 57 & 9 & 9, 57 \\
  6 & 3 & 3, 6 & 58 & 16 & 16, 58 \\
  10 & 4 & 4, 10 & 62 & 17 & 17, 62 \\
  14 & 5 & 5, 14 & 65 & 65 & 65\\
  15 & 15 & 15 & 66 & 18 & 18, 21, 33, 66\\
  21 & 5 & 5, 21 & 69 & 69 & 69\\
22 & 7 & 7, 22 & 70 & 19 &  19, 25, 28, 70\\
26 & 8 & 8, 26 & 74 & 20 & 20, 74  \\
30 & 9 & 9, 12, 15, 30 & 77 & 77 & 77\\
33 & 33 & 33 & 78 & 8 &  8, 14, 21, 24, 39, 78 \\
34 & 10 & 10, 34 & 82 & 22 & 22, 82 \\
35 & 35 & 35 & 85 & 85 & 85\\
38 & 11 & 11, 38 & 86 & 23 & 23, 86\\
39 & 7 & 7, 39 & 87 & 87 & 87\\
42 & 7 & 7, 10, 12, 15, 21, 42 & 91 & 91 & 91 \\
46 & 13 & 13, 46 & 93 & 13 & 13, 93\\
51 & 51 & 51 & 94 & 25 & 25, 94\\
55 & 7 & 7, 55 & 95 & 95 & 95\\

 \hline
\end{tabular}
\end{table}

{\bf Example: }Let us, in order to better understand what the formula
we have derived means and how
it can be used, by hand compute the smallest class number that a group
of order $n=2 \cdot 3 \cdot 5 \cdot 7 = 210$ can have. We have to
consider the different factorisations $n=de$ and in each case compute the
smallest class number. It is not hard to see that the values in the
$\alpha$-sequence are minimal when the values in the $d$-sequence are
maximal and that this occurs when $d_i=(d,p_i)$. (See the beginning of
section \ref{formulasec} for the definition and basic properties of the
$d$- and $\alpha$-sequences.) This gives us the
smallest class number for a given $d$. Moreover, not all
factorisations $n=de$ need to be considered. Recall that the group
for which our formula gives the class number has a presentation
(\ref{Gpres}) where $t$ is an integer which is of order $d$ modulo
$e$ and that this implies that $d \mid \phi(e)$ by Remark
\ref{decond}.

In our case the
divisors of $n$ are $d=1,2,3,5,7,6,10,14,15,21,35,30,42,70,105,210$.
Out of those only $d=1,2,3,6$ satisfy $(d,\phi(e))=d$. The smallest
class number of a group corresponding to a certain $d$ can now be
computed by determining first $d_i=(d,p_i)$ (where $e=\prod_{i=1}^m p_i$) and
then $\alpha_j=\prod_{d_i|j}p_i$ which is substituted into
(\ref{formula}). We may omit the computations for $d=1$ as $G$ is this
case is cyclic and hence has class number $|G|=210$.

For $d=2$ we have $e=p_1p_2p_3=3 \cdot 5 \cdot 7$, and hence
$d_1=(d,p_1-1)=(2,2)=2, d_2=(d,p_2-1)=(2,4)=2$ and
$d_3=(d,p_3)=(2,6)=2$ so ${\bf d}=(2,2,2)$. From this we get $\alpha_1
= \prod_{d_i|1}p_i=1$ and $\alpha_2=\prod_{d_i|2}p_i=3 \cdot 5 \cdot
7=105$ so ${\bf \alpha}=(1,105)$. Now by formula (\ref{formula}) the
class number is the dot product of ${\bf \alpha}$ and
$2(1-\frac{1}{2^2},\frac{1}{2^2})=\frac{1}{2}(3,1)$ so
$k(G)=\frac{1}{2}(3,1) \cdot (1,105)= 54$

A similar computation for $d=3$ gives ${\bf d}=(1,1,3)$, ${\bf
  \alpha}=(10,70)$ and
$k(G)=3(10,70) \cdot (1-\frac{1}{3^2},\frac{1}{3^2})=\frac{1}{3}(10,70) \cdot
(8,1) = 50$.

Finally, for $d=6$ we have ${\bf d}=(2,6)$, ${\bf \alpha}=(1,5,1,35)$ and
$k(G)=6(1,5,1,35)(1-\frac{1}{2^2}-\frac{1}{3^2}+\frac{1}{2^23^2},\frac{1}{2^2}(1-\frac{1}{3^2}),\frac{1}{3^2}(1-\frac{1}{2^2}),\frac{1}{2^23^2})=\frac{1}{6}(1,5,1,35)\cdot
(24,8,3,1)=17$.

We have now computed the smallest class number for each choice of $d$
and obtained $k(G)=210,54,50,17$ for $d=1,2,3,6$ respectively. We
conclude that the smallest possible class number for a group of order
$210$ is $17$ and that this occurs for the group with presentation

$$ \langle x,y|x^6=id,y^35=id,x^{-1}yx=y^t \rangle $$

where $t$ is a number of order $6$ modulo $35$, order $2$ modulo $5$
and order $6$ modulo $7$. There are two possible values, $t=19,24$,
which give isomorphic groups.

$\blacklozenge$

\vspace{4mm}

As the above example illustrates the problem when wanting to find the
smallest class number that can occur for groups of a given order $n$ is
that we do not know in advance which factorisation $n=de$ will give
the smallest class number. Thus we would like to know how to choose
$d$ to make the class number minimal. As long as the number of factors
in $n$ is small one can analyse all possibilities.

Let us see what happens if $n=pqr$ where $p < q < r$ are primes. The
class number is determined from our formula once we have computed the
$\alpha$-sequence, and this sequence is in turn determined by the
${\bf }d$-sequence. As we have noted we only have to consider the
choice $d_i=(d,\phi(p_i))$, $p_i$ being the $i$th factor in $e$, if
searching for the smallest class number a group with presentation
(\ref{Gpres}) can have. This means that all information we need about
$p$, $q$ and $r$ is something that allows us to compute the
$d_i$. This can be done if we know exactly which of the three
conditions below are satisfied
\begin{itemize}
  \item {\bf A:} $p|q-1$
  \item {\bf B:} $p|r-1$
  \item {\bf C:} $q|r-1$
\end{itemize}

Considering all combinations of factorisations $n=de$ and conditions
we can compute all class numbers for groups of order $n=pqr$ and then
for each condition situation see which choice of factorisation
turns out to be optimal. The resulting class numbers are shown in
table \ref{pqrtabell}.

\noindent
\begin{table}
  \caption{Class numbers for groups of order $pqr$. The capital
    letters in the leftmost column describes which of the conditions
    ${\bf A:} p|q-1$,  ${\bf B:} p|r-1$ and ${\bf C:} q|r-1$ are
    satisfied. } \label{pqrtabell}
   \vspace{2mm}
\begin{tabular}{|l||l|l|l|l|l|l|l|l|}
  \hline
 & $d=1$ & $d=p$ & $d=q$ & $d=r$ & $d=pq$ & $d=pr$ & $d=qr$ & $d=pqr$ \\
 \hline \hline
 - & ${\bf pqr}$ & - & - & - & - & - & - & - \\
 A & $pqr$ & ${\bf pr+\frac{r(q-1)}{p}}$ & - & - & - & - & - & - \\
 B & $pqr$ & ${\bf pq+\frac{q(r-1)}{p}}$ & - & - & - & - & - & - \\
 C & $pqr$ & - & ${\bf pq+\frac{p(r-1)}{q}}$ & - & - & - & - & -  \\
 AB & $pqr$ & ${\bf p+\frac{qr-1}{p}}$ & - & - & - & - & - & - \\
 AC & $pqr$ & $ pr+\frac{r(q-1)}{p}$ & ${\bf pq+\frac{p(r-1)}{q}}$ & - &
 - & - & - & - \\
 BC & $pqr$ & $ pq+\frac{q(r-1)}{p}$ & $pq+\frac{p(r-1)}{q}$ & - & ${\bf
   pq+\frac{r-1}{pq}}$ & - & - & - \\
 ABC & $pqr$ & $ p+\frac{qr-1}{p}$ &  $pq+\frac{p(r-1)}{q}$ & - &  ${\bf
   pq+\frac{r-1}{pq}}$ & - & - & - \\
 \hline
\end{tabular}
\end{table}

As certain combinations of conditions and
factorisations do not satisfy $(d,\phi(e))=d$ there are a number of
blank spaces. The entries in boldface are the smallest in
their row. From this table it is evident that when $n$ is a product of
three primes the choice of $d$ that gives the smallest class number is
the largest $d$ satisfying $(d,\phi(n/d))=d$. It is easy to verify
that this also holds when $n$ has less than three prime factors. One
might then conjecture this to be true for all squarefree
$n$. Unfortunately this is not the case, and a smallest contradiction
is obtained for $n=930 = 2 \cdot 3 \cdot 5 \cdot 31$, a product of
four primes. The smallest class number of a group of order $930$
equals $24$ and is obtained for $d=10$. The largest $d$ with
$(d,\phi(n/d))=d$ is in this case $d=30$, but this choice of $d$
produces larger class numbers than $24$, the smallest one being $31$. A complete
analysis of the $d$-values which give the smallest class numbers when
$n$ is a product of four (or more) factors might make us able to
conjecture a general rule for the optimal choice of $d$. However there
are 64 cases to consider, so it is advisable to computerise the
process. We leave this problem open for now.

\section{Ideas for future work}

As mentioned in the previous paragraph it would be nice to be able to
find the optimal choice of $d$ as a function of $n$ as this would give
a direct formula for the smallest class number among those of groups
of a certain squarefree order.

Another view of the problem of finding small class numbers is to start
out with a certain class number and then try to construct the largest
group having this class number. (Or at least determine the order of
this group.) As mentioned in the introduction general attempts in this direction exist. Confining ourselves to groups of squarefree order we can solve
the problem as follows. We already have an algorithm to compute all
class numbers a group of order $n$ can have. By Theorem
\ref{upperbound} we know that
all groups with $k$ conjugacy classes are of order at most
$k^3$. Hence, computing the class numbers of all groups of orders up
$k^3$ and checking which of them has class numbers exactly equal to
$k$ we obtain all groups with $k$ conjugacy classes. We have performed
such a computation for $k$ up to $100$ and the results are provided in
tables \ref{tabellallgroups} and \ref{tabelllargest}.

\noindent
\begin{table}
  \caption{Squarefree group orders for which there exists a group with
    given class number.} \label{tabellallgroups}
  \vspace{2mm}
\begin{tabular}{|l|l|}
  \hline
  $k(G)$ & All squarefree $|G|$ \\
  \hline \hline
  1 & 1 \\
  2 & 2 \\
  3 & 3, 6 \\
  4 & 10\\
  5 & 5, 14, 21\\
  6 & 6 \\
7 & 7, 22, 39, 42, 55 \\
8 & 26, 78  \\
9 &  30, 57, 114 \\
10 & 10, 34, 42 \\
11 &  11, 38, 110, 155, 186, 203 \\
12 & 30, 42, 222 \\
13 &  13, 46, 93, 205, 253, 258, 301, 310 \\
14 &  14, 78, 110, 410 \\
15 &   15, 30, 42, 111 \\
16 & 58, 366, 406, 610 \\
17 & 17, 62, 129, 210, 305, 402, 465, 497, 602, 689, 710, 737 \\
18 & 66, 114, 330, 438 \\
19 & 19, 70, 355, 474, 915, 979, 994, 1027 \\
20 & 74, 210, 1010 \\
21 &  21, 42, 66, 78, 165, 330, 546 \\
22 &  22, 82, 310, 390, 406, 582, 1582 \\
23 &  23, 86, 183, 506, 618, 791, 903, 1310, 1703, 1751, 1778\\
24 &  78, 546, 654, 930 \\
25 &  70, 94, 105, 201, 505, 889, 1081, 1474, 1510, 2041, 2265, 2329 \\
26 &  26, 186, 410, 506, 602, 1958 \\
27 &  102, 219, 546, 570, 762, 1218, 1230, 2667, 2715 \\
28 & 70, 106, 390, 798, 1378, 1810, 2758\\
29 &  29, 110, 237, 546, 834, 1910, 2054, 2189, 2954, 3165, 3197, 3629 \\
30 & 30, 102, 114, 222, 1806 \\
31 &  31, 118, 655, 798, 906, 930, 1711, 2110, 3346, 3406, 3615, 4063, 4351, 4378,
4431\\
32 &  122, 462, 770, 942, 1830, 4082 \\
33 &  33, 66, 114, 273, 330, 465, 609, 930, 978, 1218, 1830, 4065 \\
34 &  34, 130, 258, 610, 930, 994, 1378, 1474, 2410, 3934\\

 \hline
\end{tabular}
\end{table}

\noindent
\begin{table}
    \caption{Largest squarefree group order for which there exists a group with
    given class number.} \label{tabelllargest}
   \vspace{2mm}
\begin{tabular}{|l|l|l|l|}
  \hline
  $k(G)$ & Largest squarefree $|G|$ &  $k(G)$ & Largest squarefree $|G|$\\
  \hline \hline
35 & 6371 & 68 & 46598 \\
36 & 5430 & 69 & 46869 \\
37 & 7282 & 70 & 46246 \\
38 & 8138 & 71 & 53063 \\
39 & 8130 & 72 & 26202 \\
40 & 8734 & 73 & 57481 \\
41 & 10121 & 74 & 59294 \\
42 & 6510 & 75 & 26106 \\
43 & 11563 & 76 & 63526 \\
44 & 12630 & 77 & 67677 \\
45 & 13101 & 78 & 31314 \\
46 & 13906 & 79 & 72879 \\
47 & 15279 & 80 & 74066 \\
48 & 16230 & 81 & 78081 \\
49 & 17185 & 82 & 76222 \\
50 & 18030 & 83 & 84531 \\
51 & 18930 & 84 & 48630 \\
52 & 20842 & 85 & 91029 \\
53 & 22085 & 86 & 92318 \\
54 & 10218 & 87 & 73047 \\
55 & 24586 & 88 & 100978 \\
56 & 21206 & 89 & 104249 \\
57 & 26502 & 90 & 84714 \\
58 & 28918 & 91 & 111691 \\
59 & 30299 & 92 & 115382 \\
60 & 31794 & 93 & 117678 \\
61 & 33661 & 94 & 121162 \\
62 & 32402 & 95 & 126815 \\
63 & 37086 & 96 & 95298 \\
64 & 38134 & 97 & 135265 \\
65 & 40721 & 98 & 139458 \\
66 & 42378 & 99 & 143814 \\
67 & 10402 & 100 & 146134 \\

  \hline
\end{tabular}
\end{table}

 It is likely that
large group orders occurring in this table display certain properties
of their prime factorisations, since this is what influences the ${\bf
  d}$-sequence determining the class number. An idea for future work
is to try to pinpoint those properties.

In the bigger picture this work on groups of squarefree order is one
contribution to the greater challenge of understanding class numbers
of groups in general. Obviously it would be interesting to try to
generalise our results. The ultimate goal would be to see if a
modified version can apply to groups in general, but perhaps
more within reach if we can cover groups with presentations similar to (\ref{Gpres}).

\section* {Acknowledgments}

I would like to express my gratitude to Professor John McKay for inviting me to Concordia University and
for introducing me to the fascinating problem of class numbers of
groups. Many thanks also to Olof Barr for valuable suggestions and
comments during my work with this article.

\section*{Appendix}

\begin{verbatim}

// Computation of all possible class numbers for given group order
//    n (n squarefree). The results are stored in Classtable. The
//    smallest class number for each n is stored in MinClasstable


Classtable:=[];
Lowerlimit:=2; Upperlimit:=100;
for n in [Lowerlimit..Upperlimit] do
  if IsSquarefree(n) then

    Classnrs:=[];

for d in Divisors(n) do
  if (GreatestCommonDivisor(d,EulerPhi(n div d))
  eq d) then
e:=n div d;

facte:=Factorisation(e);
m:=#facte;
MaxOrders:=[];
for i in [1..m] do
  MaxOrders:=Append(MaxOrders,GreatestCommonDivisor(d,facte[i][1]-1));
end for;
PossOrders:=[];
for i in [1..m] do
 PossOrders:=Append(PossOrders,Divisors(MaxOrders[i]));
end for;
  OrderSeqs:=[];
  for s in [1..#PossOrders[1]] do
  OrderSeqs:=Append(OrderSeqs,[PossOrders[1][s]]);
  end for;

  for j in [2..m] do
    NewOrderSeqs:=[];
    for k in [1..#PossOrders[j]] do
    for O in OrderSeqs do
      NewOrderSeqs:=Append(NewOrderSeqs,Append(O,PossOrders[j][k]));
    end for;
    end for;
    OrderSeqs:=NewOrderSeqs;
  end for;


alphaseqs:=[];

for k in [1..#OrderSeqs] do
alphaseq:=[];
for j in [1..d] do
  alpha:=1;
  for t in [1..#OrderSeqs[k]] do
    if ((j mod (OrderSeqs[k][t])) eq 0) then
      alpha:=alpha*facte[t][1];
    end if;
  end for;
  alphaseq:=Append(alphaseq,alpha);
end for;
alphaseqs:=Append(alphaseqs,alphaseq);
end for;

for alphaseq in alphaseqs do
  sum:=0;
  for j in Divisors(d) do
    insum:=1;
    for s in Divisors(d div j) do
      r:=#Factorisation(s);
      if not s eq 1 then
        insum:=insum+(-1)^r/s^2;
      end if;
    end for;
    insum:=insum*alphaseq[j]/j^2;


    sum:=sum+insum;
  end for;
  sum:=d*sum;
Classnrs:=Append(Classnrs,sum);
end for;
end if;
end for;

// Classnrs contains  all class numbers for
//groups of order n computed by our formula.

Classtable[n]:=SequenceToSet(Classnrs);
end if;
end for;

\end{verbatim}
\end{document}